\documentclass[12pt]{article}

\usepackage{latexsym}
\usepackage{amsmath,amsfonts}
\usepackage[dvips]{epsfig}
\usepackage[margin=1in]{geometry}

\newcommand{\p}{\partial}
\newcommand{\e}{\epsilon}
\newcommand{\n}{\nabla}
\newcommand{\om}{\int_\Omega}

\newcommand{\dv}{{\rm div}}

\newcommand{\be}{\begin{equation}}
\newcommand{\ee}{\end{equation}}

\newtheorem{thm}{Theorem}
\newtheorem{defn}{Definition}
\newtheorem{lem}{Lemma}

\newtheorem{rmk}{Remark}

\numberwithin{equation}{section}
\numberwithin{thm}{section}
\numberwithin{lem}{section}
\numberwithin{defn}{section}
\numberwithin{prop}{section}
\numberwithin{rmk}{section}

\begin{document}




\title{Incompressible limit of strong solutions to 3-D Navier-Stokes equations with Navier's slip boundary condition for all time\footnote{This work was initiated when the author Ou visited The Institute of Mathematical Science, The Chinese University of Hong Kong in 2012. He would like to thank Prof. Zhouping Xin for the invitation and helpful discussions.}}


\author{Yaobin Ou$^{1,3}$, Dandan Ren$^2$\\
\small{${}^1$School of Information, Renmin University of China, Beijing 100872, P.R. China}
\\
\small{${}^2$School of Mathematical Sciences, University of Electronic Science and Technology of China,}\\
\small{ Chengdu, Sichuan 611731, P.R. China}
\\
\small${}^3${The Institute of Mathematical Sciences, The Chinese University of Hong Kong,}\\
\small{Shatin, N.T., Hong Kong}\\
\small{Email: {ou.yaobin@gmail.com, rendanlengyan@163.com}}
}
\maketitle
\begin{abstract}
This paper studies the incompressible limit of global strong solutions to the three-dimensional
compressible Navier-Stokes equations associated with Navier's slip boundary condition, provided that the time derivatives, up to first order, of solutions are bounded initially. The main idea is to derive a differential inequality with decay, so that the estimates are bounded uniformly both in the Mach number $\epsilon\in (0,\epsilon_0]$ (for some $\epsilon_0>0$) and the time $t\in[0,+\infty)$.
\end{abstract}




\section{Introduction}
The motions of  highly subsonic viscous fluids in a bounded domain $\Omega\subset\mathbb{R}^3$ are described by the following non-dimensionalized Navier-Stokes equations:
\begin{eqnarray}
 \label{NS1}
&&\rho_t+\textrm{div}(\rho u)=0,\\
 &&(\rho u)_t+\textrm{div}(\rho u\otimes
u)-\textrm{div}\mathbb{S}+\frac{1}{\e^2}\nabla p=0,\label{NS2}
\end{eqnarray}
where the first equation represents the conservation of mass and the second one denotes the conservation of momentum.   The unknowns $\rho$, $u$ and $p$ are the density, the velocity and the pressure, respectively. And the matrix
$\mathbb{S}\equiv\mathbb{S}(u)=2\mu D(u)+\lambda{\rm div} \textrm{u}I_3$ is the viscous stress tensor for newtonian fluids, where $D(u)=(\nabla u+\nabla u^t)/2$. The constant $\epsilon\in (0,1]$ is the Mach number of the highly subsonic fluids. The constants $\mu,\;\lambda$ are viscosity coefficients with $\mu>0$, $\mu+3\lambda/2\ge0$. In this paper, we suppose that the pressure
$p=p(\rho)$ is a $C^3$ function satisfying that $p'(\rho)>0$ for
$\rho>0$ .

In the physical viewpoint, the motions of highly subsonic
compressible fluids would behave similarly to the incompressible ones (see
\cite{LL}). Formally, as the Mach number $\e$ tends to zero, the
solutions to \eqref{NS1}-\eqref{NS2} will converge to the solution
$(u,\pi)$ of the incompressible Navier-Stokes equations, namely
\begin{eqnarray*}
&& u_t+u\cdot\nabla u-\mu\Delta u+\nabla\pi=0,
\\
&&{\rm div}u=0.
\end{eqnarray*}
It is known as the incompressible limit, which is one of the fundamental hydrodynamic limits. However, the rigorous justification of the limit poses challenging problems mathematically since singular phenomena usually occur in this process. To be precise, both the uniform estimates in Mach number and
the convergence to the incompressible model are usually difficult to
obtain. In the following, we restrict the discussion in the isentropic regime only.

The general framework for studying the incompressible limit for local strong or smooth solutions was established by S. Klainerman and A. Majda in \cite{KM1,KM2}. In these works, they proved the incompressible limit of local smooth solutions to the Navier-Stokes equations (or the Euler equations) with \lq\lq
well-prepared" - some smallness assumption on the divergence of
initial velocity - initial data, in $\mathbb{R}^n$ or $\mathbb T^n$.
Indeed, by analyzing the rescaled linear group generated by the penalty operator of order $\epsilon^{-1}$ (see \cite{SC1,UK} for instance), the incompressible limit can also be verified for the cases of general data that the velocity of incompressible fluid is just the limit of Leray projection for the velocities in compressible fluids.
This method also applies to global weak solution of the isentropic Navier-Stokes
equations with general initial data and various boundary conditions
\cite{DG,DGLM,LM}. Especially, P.-L. Lions and N. Masmoudi \cite{LM}  studied
 the incompressible limit for the weak solutions to the Navier-Stokes equations with a slip boundary condition, that is, on the boundary $\partial\Omega$ of $\Omega\subset \mathbb{R}^n$,
 \begin{eqnarray}
 && u\cdot n=0,\; {\rm curl} u=0\quad{\rm for}\quad n=2,\quad {\rm or} \label{bc2d}\\
 && u\cdot n=0,\; n\times {\rm curl} u=0 \quad{\rm for}\quad n=3,
 \end{eqnarray}
 where ${\rm curl} u=(\partial_2 u, -\partial_1 u)^t $ for $n=2$ and ${\rm curl} u=(\partial_2u_3-\partial_3u_2,\partial_3u_1-\partial_1u_3,\partial_1u_2-\partial_2u_1)^t$ for $n=3$. Recently, D.~Donatelli, E.~Feireisl, A. Novotn\'y, etc. have also obtained a series of important progresses on incompressible limits of weak solutions to compressible Navier-Stokes equations associated with slip boundary conditions (see \cite{DFN,Fe}, for instance).
 For other interesting results  on the incompressible limit in a finite time interval, which may be independent of the initial data, for isentropic fluids, the reader may refer to
\cite{DAN1,HW,JJL,Lin,Ma1,SEC1,SEC3} and many others.

Although numerous significant progresses on incompressible limit had been achieved during the last four decades, only a few results were concerned with  global strong or classical solutions for the time $t\in [0,+\infty)$. In this situation, one needs to show the uniform estimates with respect to both $\epsilon\in (0,\epsilon_0]$ (for some small constant $\epsilon_0>0$) and $t\in [0,+\infty)$. Thus additional difficulties arise.  D. Hoff \cite{Hoff} verified the incompressible limit for the global solutions in $\mathbb{R}^3\times[0,+\infty)$ with general initial data, provided that the background solution to the incompressible Navier-Stokes equations is sufficiently smooth.
For regular solutions with no-slip boundary conditions , i.e., $u|_{\partial\Omega}=0,$ where $\Omega\subset \mathbb{R}^3$ is a bounded domain, and slightly compressible initial data, H. Bessaih \cite{Be} established the uniform estimates both in the Mach number and $t\in[0,+\infty)$, and showed the strong
convergence to the solution of incompressible Navier-Stokes equations. In \cite{OU1}, the author studied the incompressible limit of regular
solutions to the compressible Navier-Stokes equations \eqref{NS1}-\eqref{NS2} with slightly compressible initial data in a 2-D bounded domain with the boundary
condition in \eqref{bc2d}.

The aim of this paper is to extend the result in \cite{OU1} to three spatial dimensions, that is, to study the incompressible limit of global strong
solutions to the 3-D compressible Navier-Stokes equations \eqref{NS1}-\eqref{NS2} with  Navier's slip boundary condition
\be\label{Navier}
u\cdot n=0,\quad \tau\cdot\mathbb{S}(u)\cdot n+\alpha
u\cdot\tau=0\quad {\rm on} \; \p\Omega,
\ee
where $n,\tau$ are the unit outer normal and  tangential vector to the boundary, respectively.
This is a non-trivial generalization since on the boundary of a 3-D bounded domain, the information on the normal component of the vorticity curl$u$ is unavailable (see Lemma \ref{lem25} for instance), thus the classical regularity theory dosn't apply.

At the same time, this paper also generalizes the result in \cite{Be} in the sense that all the second-order spatial derivatives are uniformly bounded with respect to the Mach number and the time.
Moreover, it is worthy to note that the result in this paper can not be covered by the ones in \cite{DFN,Fe} since our estimates are uniformly bounded for all the time in $[0,+\infty)$, instead of a fixed finite interval. Furthermore, the method in the current paper can also simplify the proof of  local existence results in \cite{Za}.

To simplify the proof, we convert the equations into the anti-symmetric form by setting $\rho=1+ \epsilon \sigma$. Then the Navier-Stokes equations \eqref{NS1}-\eqref{NS2} are equivalent to
\begin{equation} \label{NS3}
\sigma_t+\textrm{div}(\sigma u)+\frac{1}{\epsilon} \textrm{div}u=0,
\end{equation}
\begin{equation} \label{NS4}
\rho(u_t+u\cdot\n u)+\frac{1}{\epsilon}p'(1+\e \sigma)\n \sigma=2\mu\textrm{div}(D(u))+\lambda\n{\rm div} \textrm{u}.
\end{equation}
For the new unknowns $(\sigma, u)$, we impose the following initial condition
\be\label{NS5}
(\sigma,u)|_{t=0}=(\sigma_0,u_0)(x),\;x\in\Omega
\ee
and the slip boundary condition
\be\label{NS6}
u\cdot n=0,\quad \tau\cdot D(u)\cdot n+\alpha
u\cdot\tau=0\quad {\rm on} \; \p\Omega,
\ee
which is equivalent to \eqref{Navier}. One may refer to
\cite{Za}, for instance, for the description  and the background on this boundary condition.

First,  the  local existence of the solution $(\sigma,u)$ to the problem \eqref{NS3}-\eqref{NS6} is indeed established by W.M. Zajaczkowski \cite{Za} in the framework of \cite{MN2,V}.

{\thm {\rm (Local existence)}\label{th11} Let $\epsilon\in(0,1]$ be a fixed constant and $\Omega\subset\mathbb{R}^3$ be a simply connected, bounded domain with smooth boundary $\partial\Omega$. Suppose that the initial datum $(\sigma_0,u_0)$ satisfies the following conditions,
$$(\sigma_0,u_0) \in H^2(\Omega), \quad (\sigma_t(0),u_t(0)) \in H^1(\Omega), \quad (\sigma_{tt}(0),u_{tt}(0)) \in L^2(\Omega),$$
with $\om \sigma_0 dx = 0$ and $1+\epsilon\sigma_0 \ge m$ for some positive constant $m$. Assume the following compatibility conditions are satisfied:
$$
\partial_t^i u(0)\cdot n = \tau\cdot\mathbb{S}(\partial_t^i u(0))\cdot n+\alpha\tau\cdot
\partial_t^i u(0) = 0 \quad {\rm on} \; \p\Omega, \; i=0,1,2.
$$
Then there exists a positive constant $T = T(\sigma_0,u_0,m,\epsilon)$ such that the initial-boundary problem \eqref{NS3}-\eqref{NS6} admits a unique solution $(\sigma,u)$  satisfying that $1+\epsilon\sigma_0>0$ in $\Omega\times(0,T)$, and
\begin{equation*}
\begin{split}
&\sigma\in C([0,T],H^2), \quad u\in C([0,T],H_0^1\cap H^2)\cap L^2(0,T;H^3), \\
&\sigma_t\in C([0,T],H^1), \quad u_t\in C([0,T],H^1_0)\cap L^2(0,T;H^2),\\
& \sigma_{tt}\in C([0,T],L^2), \quad u_{tt}\in C([0,T],L^2)\cap
L^2(0,T;H^1 ).
\end{split}
\end{equation*}
}
\noindent{\rmk To simplify the statement, we use the notation \lq\lq$u_t
(0)$" to signify the quantity $u_t |_{t=0}:=-u_0 \cdot\nabla
u_0  + (-p'(1+\e \sigma_0)\n \sigma/\epsilon +\mu\Delta u_0  +\nu\nabla\dv\,
u_0 )/(1+\e \sigma_0) $ obtained from the equation \eqref{NS4}. And the notation \lq\lq$\partial_t^i u(0)$"  is given  by differentiating \eqref{NS4} $i-1$ times with respect to $t$ and then letting $t=0$. The same rule applies to the notations $\partial_t^i \sigma(0),\;\partial_t^i \rho(0)$, etc.}
\vskip .2cm
The purpose of this paper is to prove the following uniform estimates with respect to $\epsilon\in (0,
 \epsilon_0]$ (for some $0<\epsilon_0\le 1$) and $t\in (0,+\infty)$, thus apply Theorem \ref{th11} to obtain the global existence theorem and the corresponding incompressible
limits. In order to state the theorem precisely, we introduce the following notation
{\defn
$$\phi^\e(t) := \max_{0\le s\le t}(\|(\sigma^\e,u^\e)\|_{H^2}+\|(\sigma_t^\e,u_t^\e)\|_{H^1} +\e \| (u_{tt}^\e,\sigma_{tt}^\e)\|_{L^2})(s).$$
}

Then the main results of this paper is stated as follows.
{\thm {\rm (Global-in-time existence and incompressible limit)}. Let all the assumptions
in Theorem \ref{th11} be satisfied and $\epsilon\in(0,\overline{\epsilon}]$ for some sufficiently small constant $\overline{\epsilon}\in (0,1]$. Moreover, we assume that
$$
\phi^\e(0) \le \theta,
$$
for some sufficiently small positive constant $\theta$. Then there exists a unique solution $(\sigma^\e,u^\e)$ to the initial-boundary value problem  \eqref{NS3}-\eqref{NS6} in $\Omega\times\mathbb{R}^+$, such that
\begin{equation*}
\begin{split}
&\sigma^\e\in C(\mathbb{R}^+,H^2), \quad u^\e\in C(\mathbb{R}^+,H_0^1\cap H^2)\cap L^2(\mathbb{R}^+;H^3), \\
&\sigma_t^\e\in C(\mathbb{R}^+,H^1), \quad u_t^\e\in C(\mathbb{R}^+,H^1_0)\cap L^2(\mathbb{R}^+;H^2),\\
& \sigma_{tt}^\e\in C(\mathbb{R}^+,L^2), \quad u_{tt}^\e\in C(\mathbb{R}^+,L^2)\cap
L^2(\mathbb{R}^+;H^1 ).
\end{split}
\end{equation*}
where $\mathbb{R}^+ = [0,+\infty)$. Furthermore, the following uniform estimate in $\epsilon\in(0,\overline{\epsilon}]$ holds:
\be\label{phit}
\phi^\e(t) \le C\theta, \quad \forall t \in\mathbb{R}^+,
\ee
for some positive constant $C$. Thus $u^\e\rightarrow v$ strongly in $C(\mathbb{R}_{loc}^+;H^s)$ as $\epsilon\rightarrow 0$ for any $0 \le s<2$. And there exists a function $P(x,t)$, such that $(v,P)$ is the unique solution of the following initial-boundary value problem of the incompressible Navier-Stokes equations:
\begin{eqnarray*}
&& {\rm div}v = 0,\; v_t + v\cdot\n v + \n P = \mu\Delta v,\quad {\rm in}\quad \Omega\times (0,+\infty),\\
&& v\cdot n = \tau\cdot\mathbb{S}(v)\cdot n+\alpha\tau\cdot v = 0 \quad {\rm on}\; \partial\Omega,\\
&& v|_{t=0} = v_0(x), \; x\in\Omega,
\end{eqnarray*}
where $\|u_0^\e-v_0\|_{H^2}+\|u_0^\e-v_0\|_{H^1}+\|\e(u_0^\e-v_0)\|_{L^2}\rightarrow0$ as $\epsilon\rightarrow0$.
}
\vskip .2cm
\noindent{\bf Proof.} This theorem was shown by Lemma \ref{lem58} and the same arguments as in \cite{Be,OU1}.

{\rmk Although the time derivatives up to second order are estimated, however, only the derivatives up to first order are required to be bounded initially.}
\vskip .2cm
The rest of this paper is organized as follows. In Section 2, we present some lemmas which will be used in estimating the Sobolev norms in a bounded domain and dealing with the slip boundary condition. In Section 3, we show the uniform-in-$\epsilon$ estimates  by deriving a differential inequality with certain decay property.  We first show the $L^2$ estimate of the solutions, next the low-order spatial, temporal or mixed  derivatives, and then the high-order derivatives.
 The strategy for estimating derivatives is to treat the vorticity and the divergence
of velocity respectively, based on the decomposition $\triangle=\n{\rm div} -{\rm curl} {\rm curl}$ and the slip boundary condition.  Moreover, to overcome the difficulty in estimating the vorticity, due to the loss of information on the normal component, we  take the advantage of the isothermal coordinates to estimate it in local regions near the boundary.
 By combining carefully all the spatial-temporal estimates, we obtain the uniform estimate with respect to both $\epsilon\in (0,\bar\epsilon] \;(0<\bar\e\le 1)$  and $t\in [0,+\infty)$.
.

\section{Preliminaries}
Throughout this paper, we will use the following lemmas from time to time.
{\lem {\rm(See \cite{BB})}\label{p2}. Let $\Omega$ be a bounded domain in $\mathbb{R}^N$ with smooth boundary $\partial \Omega$ and outward normal $n$. Then there exists a constant $C>0$ independent of $u$, such that
$$
\|u\|_{H^s(\Omega)} \le C(\|{\rm div} u\|_{H^{s-1}(\Omega)}+\|{\rm curl} u\|_{H^{s-1}(\Omega)}+\|u\cdot n\|_{H^{s-\frac{1}{2}}(\partial\Omega)}+\|u\|_{H^{s-1}(\Omega)}),
$$
for any $u \in H^s(\Omega)^N$.
}

{\lem\label{lem26}{\rm(See \cite{BB}).} Assume $f \in C([0,T];W^{k,p}(\Omega,\mathbb{R}^N))$ with
$$k > \frac{N}{p}+1 \quad and \quad 1 \le p \le +\infty.$$
Then the problem
$$\frac{du}{dt}(x,t) = f(u(x,t),t), \quad u(x,0) = x$$
has a solution $u \in C^1([0,T];D^{k,p}(\Omega))$, where
$$
D^{k,p}(\Omega) = \{\eta \in W^{k,p}(\Omega) \;| \;\eta \;is \;a \;bijective \;from \; \overline{\Omega} \;onto \;\overline{\Omega},\; \eta^{-1} \in W^{k,p}(\Omega)\}.$$
}
{\lem\label{lem27}{\rm(See \cite{BB}).} Let $k \geq 2$ be an integer, and let $1 \le p \le q \le +\infty$ be such that $p < +\infty$ and $k > \frac{N}{p}+1$. Let $f \in W^{k,p}(\Omega)$, then the mapping $g \longmapsto g \circ f$ is continuous from $D^{k,p}(\Omega))$ into $W^{k,p}(\Omega)$.
}

{\lem {\rm(See \cite{XX})}.\label{p3} Let $\Omega$ be a bounded domain in $\mathbb{R}^N$ with smooth boundary $\partial U$ and outward normal $n$. Then there exists a constant $C>0$ independent of $u$, such that
$$
\|u\|_{H^s(\Omega)} \le C(\|{\rm div} u\|_{H^{s-1}(\Omega)}+\|{\rm curl} u\|_{H^{s-1}(\Omega)}+\|u\times n\|_{H^{s-\frac{1}{2}}(\partial\Omega)}+\|u\|_{H^{s-1}(\Omega)}),
$$
for all $u \in H^s(\Omega)^N$.
}

{\lem\label{lem1} {\rm(See \cite{DV})}. Let $\Omega\subset \mathbb{R}^3$ be a
 open bounded domain with $C^2$ boundary $\p\Omega$. Moreover, we assume that $\Omega$ is simply connected and non-axisymmetric. Then for any $u\in
H^1(\Omega)$ satisfying $u\cdot n|_{\partial\Omega}=0$, one has
\begin{equation}\label{NS7}
\|u\|_{H^1(\Omega)}\le C(\|D(u)\|_{L^2(\Omega)}+\|u \|_{L^2(\p\Omega)})
\end{equation}
and
\begin{equation}\label{NS71}
\|\n u\|_{L^2(\Omega)}\le C(\|{\rm div} u\|_{L^2(\Omega)}+\|{\rm curl} u\|_{L^2(\Omega)}),
\end{equation}
}
where  $C$ is a constant independent of $u$.\\

The following lemma is a variant of Theorem 3.10 in \cite{XX2} in the case of Navier's slip boundary condition. It plays a key role in proving the vorticity estimates.
{\lem \label{lem25}. Let $\Omega\subset \mathbb{R}^3$ be a
 open bounded domain with $C^2$ boundary $\p\Omega$. If $u\in H^2(\Omega)^2$ with \eqref{NS6} being satisfied, then  \be\label{bc2}
 \tau\cdot(n\times w)  =  2\tau\cdot (\alpha u -\nabla n u ).
\ee
}

\noindent{\bf Proof.} Using the density in $\{u\in H^2(\Omega)^2| u\cdot n=0\; {\rm on}\; \partial\Omega\}$ of the velocity fields  $u\in C^\infty(\Omega)^2$ such that $u\cdot n=0$ on $\partial\Omega$, and the continuity of the trace operators, it suffices
to handle the case where $u$ is a smooth velocity field on $\bar\Omega$. Now, after extending $n(x)$ to
a tubular neighbourhood of $\partial\Omega$, we obtain
 $$\p_n u=\frac 1 2 w\times n + D(u)n$$
and
$$\p_\tau u=\frac 1 2 w\times \tau +D(u)\tau,$$
which yield
$$2(D(u)n)\cdot\tau=\p_\tau u\cdot n+\p_n u\cdot\tau$$
and
$$(n\times w)\cdot \tau=\p_\tau u\cdot n-\p_n u\cdot\tau.$$
It follows that
$$
2(D(u)n)\cdot\tau+(n\times w)\cdot \tau=2\p_\tau u\cdot n=-2u\cdot\p_\tau n=-2 \tau\cdot(\nabla n u),
$$
due to the boundary condition $u\cdot n=0$. By use of \eqref{NS6}, we easily get \eqref{bc2}.

\section{Energy estimates}

In this section, we shall derive the uniform estimates with respect to both  the time $t\in[0,+\infty)$ and the Mach
number $\epsilon\in (0,\bar\epsilon]$ for some $\bar\epsilon\in (0,1]$, which is stated as in Lemma \ref{lem58}. {We will drop the superscript $\e$ of
$\sigma^\e, u^\e, p^\e$, and so on, for the sake of simplicity.} From now
on, the positive constants $C$, $C_i$ for $i=0,1,\cdots$ below
depend only on $\Omega,$ $\mu$, $\lambda$, and $p$, but not on $T$
and $\e$. We will use $\delta$, $\eta$, and $\eta_i$ for
$i=1,2,\cdots$ to denote various small positive constants and
$C_\delta$, $C_\eta$ to denote various positive constants depending
on $\delta$ and $\eta$ respectively. For the sake of simplicity, we denote the partial derivatives $\frac{\partial}{\partial x_i}$ by $\partial_i$, $\frac{\partial^2}{\partial x_i\partial x_j}$ by $\partial_{ij}$, and so on.

 Suppose that $(\sigma,u)$ solves the initial-boundary value problem  \eqref{NS3}-\eqref{NS6} in $\Omega\times(0,T)$, for $0<T<+\infty$. In the energy estimates, we always assume that $\frac{1}{4}\leq 1+\epsilon\sigma\leq4$ in any $(x,t)\in\Omega\times(0,T)$ where $\epsilon\in(0,1]$.

 We will derive a
differential inequality in the form that,
\begin{equation*}
\frac{d}{dt}\Phi(t)+\Psi(t)\le C \Psi(t)(\Phi(t)+\Phi^2(t)),\quad \forall \,0\le t\le T,\label{decay}
\end{equation*}
where $C\ge 1$ is a constant, and $\Phi(t)$ is an equivalent norm to $\phi(t)$. Here
$\Psi(t)$ and $\Phi(t)$ are both non-negative quantities with $\Psi(t)\ge \bar C\Phi(t)$ for some constant $\bar C\in(0,1]$. The above inequality is equivalent to
$$\frac{d}{dt}\Phi(t)\le -\Psi(t)(1-C(\Phi(t)+\Phi^2(t))),\quad \forall \,0\le t\le T.$$
Thus,  if $\Phi(0)$ is small enough,  $\Phi(t)$ will be dominated by $\Phi(0)$.

\subsection{The basic estimate}

{\lem\label{z1}
There exist positive constants $C_0$ and $C_1$, such that
\begin{equation}\label{phi0}
\begin{split}
\frac{d}{dt}\Phi_0(t)+\Psi_0(t)\leq C_1 \|\sigma\|_{H^2}^2 \|\sigma\|_{H^1}^2+\epsilon\|\sigma_t\|_{L^2}^2,
\end{split}
\end{equation}
where $\Phi_0(t)=\|(\sqrt{\rho}u,\sqrt{p'(\rho)}\sigma)\|_{L^2}^2$ and $\Psi_0(t)=C_0\|u\|_{H^1}^2$.
}

\noindent {\bf Proof.} Due to the boundary conditions \eqref{NS6} and Lemma \ref{lem1}, we have
\begin{equation*}
\begin{split}
-\int_\Omega&(2\mu
\textrm{div} (D(u))+ \lambda \nabla \textrm{div}u)\cdot u dx\\
&=\om (2\mu |D(u)|^2 +\lambda({\rm div}u)^2) dx +\int_{\partial\Omega}\alpha |u|^2 dS\ge
\gamma_0 \|u\|_{H^1}^2.
\end{split}
\end{equation*}
Multiplying  \eqref{NS3} by $p'(1+\epsilon\sigma)\sigma$ and \eqref{NS4} by $u$, we get
\begin{equation*}
\begin{split}
\begin{aligned}
\frac{1}{2}\frac{d}{dt}&\|(\sqrt{\rho}u,\sqrt{p'(\rho)}\sigma)\|_{L^2}^2+\gamma_0\|u\|_{H^1}^2 \\
&\le  \frac{1}{\epsilon} \big|\int_\Omega (p'(\rho))\textrm{div}(\sigma u) dx \big| +\int_\Omega |p'(\rho)\sigma\textrm{div}(\sigma u)|dx+\frac{\epsilon}{2}\int_\Omega |p''(\rho)\sigma_t\sigma^2|dx \\
&:=A_1+A_2+A_3.
\end{aligned}
\end{split}
\end{equation*}
With the aid of the boundary condition $u\cdot n=0$, we have
\begin{eqnarray*}
A_1 &= &\big|\int_\Omega p''(\rho)\sigma\nabla\sigma\cdot u dx\big|\\
&\le & \eta\|u\|_{H^1}^2+C_\eta \|\sigma\|_{H^1}^2 \|\nabla\sigma\|_{H^1}^2.
\end{eqnarray*}
On the other hand,
$$A_2\le \eta\|u\|_{H^1}^2+C_\eta \|\sigma\|_{H^1}^2 \|\nabla\sigma\|_{H^1}^2$$
and
$$
A_3
\le \epsilon(\|\sigma_t\|_{L^2}^2+\|\sigma\|_{H^1}^4).
$$
As a result, we finish the proof of this lemma.
\hfill$\Box$

\subsection{The first-order estimate}
{\lem\label{f1}
There exists a positive constant $C_2 $ such that
$$
\begin{aligned}
\frac{d}{dt}&(\om (2\mu |D(u)|^2 +\lambda({\rm div}u)^2+\rho u_t\cdot u) dx +\int_{\partial\Omega}\alpha |u|^2 dS )+\|\sqrt{p'(\rho)}\sigma_t\|_{L^2}^2\\
&\le C_2(\|  u_t\|_{L^2}^2+ \|u\|_{H^1}^2 +  \|\sigma_t\|_{H^1}^2(\|\sigma\|_{H^2}^2+\|u\|_{H^2}^4)+\|u_t\|_{H^1}^2\|u\|_{H^1}^2).
\end{aligned}
$$
}
\textbf{Proof.} 
By differentiating  \eqref{NS4} with respect to $t$, we have
\be\label{eq32}
\begin{aligned}
(\rho u_t)_t
-& 2\mu
\textrm{div} (D(u_t))- \lambda \nabla \textrm{div}u_t+  \frac{1}{\epsilon}p'(\rho)\nabla\sigma_t\\
=& -(\epsilon\sigma_t u\cdot\n u +\rho (u_t\cdot\n u+u\cdot\n u_t)+p''(\rho)\sigma_t\n\sigma).
\end{aligned}
\ee
Then we integrate the product of \eqref{eq32} and $u$ to get
\be
\begin{aligned}\label{ddiv}
\frac{d}{dt}&\Big(\om (2\mu |D(u)|^2 +\lambda({\rm div}u)^2 + 2\rho u_t\cdot u) dx +\int_{\partial\Omega}\alpha |u|^2 dS\Big) +\frac{1}{\epsilon}\om p'(\rho)\nabla\sigma_t\cdot udx \\
&\le \|\sqrt \rho u_t\|_{L^2}^2+\eta\|u\|_{H^1}^2+C_\eta(\|\sigma_t\|_{L^2}^2(\|\sigma\|_{H^2}^2+\epsilon^2\|u\|_{H^2}^4)+\|u_t\|_{H^1}^2\|u\|_{H^1}^2).
\end{aligned}
\ee
On the other hand, we
multiply \eqref{NS3} by $p'(\rho)\sigma_t$ 
and integrate to get
\be
\begin{aligned}\label{stl}
\|\sqrt{p'(\rho)}\sigma_t\|_{L^2}^2&-\frac{1}{\epsilon}\om p'(\rho)\nabla\sigma_t\cdot u dx
\le \eta\|u\|_{H^1}^2+C_\eta\|\sigma_t\|_{H^1}^2\|\sigma\|_{H^1}^2.
\end{aligned}
\ee
Thus we summarize \eqref{ddiv} and \eqref{stl} to get the lemma.

\hfill$\Box$

From now on we may often use the following relations
$$2{\rm div}(D(u))=(\triangle u+\n\textrm{div}u),$$
$$\triangle u=\n{\rm div}u-{\rm curl} {\rm curl} u,$$
for any vector function $u=(u_1,u_2,u_3)^t$.

{\lem\label{f2}
There exists a positive constant    $C_3$ such that
$$
\begin{aligned}
\frac{d}{dt}\|\n \sigma\|_{L^2}^2&+(2\mu+\lambda)\|\sqrt{p'(\rho)^{-1}}\n{\rm div}u\|_{L^2}^2\\
& \le C_3\|u_t\|_{L^2}^2+C_3\|u\|_{H^2}^2(\|u\|_{H^1}^2+\|\sigma\|_{H^2}^2)+\eta_1\|\n^2 {\rm div}u\|_{L^2}^2+C_4\|u\|_{H^1}^2,
\end{aligned}
$$
where $\eta_1(<(2\mu+\lambda)/(8p'(4)))$ and $C_4=C_4(\eta_1)$ are to be determined  later.
}\\
\textbf{Proof.} Applying $\nabla$ to  \eqref{NS3} , we obtain
\be\label{eq37}
(\n\sigma)_t+\n^2\sigma u+\n u\n \sigma+\sigma\n{\rm div}u +{\rm div}u\n\sigma+\frac{1}{\epsilon} \n{\rm div}u=0,
\ee
then multiply the equation by $\n \sigma$ and integrate
\be
\begin{aligned}\label{sidiv}
\frac{1}{2}\frac{d}{dt}&\|\n \sigma\|_{L^2}^2+\frac{1}{\epsilon}\om \nabla\sigma\cdot\n{\rm div}u dx\\
&\le C\om(|\n u||\n \sigma|^2+|\sigma||\n{\rm div}u|| \n\sigma|)dx\\
&\le \frac{2\mu+\lambda}{8p'(4)}\|\n{\rm div}u\|_{L^2}^2+\eta\|u\|_{H^1}^2+C_\eta\|\sigma\|_{H^2}^4.
\end{aligned}
\ee
Multiplying  \eqref{NS4}  by $p'(\rho)^{-1}\n\textrm{div}u$ and integrating over $\Omega$ immediately yield
\be
\begin{aligned}\label{div}
(&2\mu+\lambda)\|\sqrt{p'(\rho)^{-1}}\n\textrm{div}u\|_{L^2}^2-\frac{1}{\epsilon}\om \n \sigma\cdot\n\textrm{div}u dx\\
&=\om p'(\rho)^{-1}(\mu \textrm{curl} \textrm{curl} u\cdot\n\textrm{div}u+\rho(u_t+u\cdot\n u)\cdot\n\textrm{div}u)dx.
\end{aligned}
\ee
Let ${\rm curl} u=w$. Then
\be\label{cudiv}
\begin{aligned}
&\om p'(\rho)^{-1} \textrm{curl} \textrm{curl} u\cdot\n\textrm{div}u dx\\
=&\om \epsilon p'(\rho)^{-2}p''(\rho)(\nabla\sigma\times w)\cdot\nabla{\rm div}udx + \int_{\partial\Omega}p'(\rho)^{-1}(n\times w)\cdot\nabla{\rm div}u dS.
\end{aligned}
\ee
Therefore, using \eqref{bc2} and the trace theorem, the boundary integral in \eqref{cudiv} can be dominated by
\begin{equation*}
\begin{aligned}
C\int_{\partial\Omega} &|(n\times w)\cdot\tau| |\nabla{\rm div}u| dS\\
&\le C\int_{\partial\Omega}2|\alpha u -2\n n u||\n\dv u|dS\\
&\le \eta \|\n {\rm div}u\|_{H^1}^2 + C_\eta\|u\|_{H^1}^2.
\end{aligned}
\end{equation*}
Thus, we can get the following inequality from \eqref{div} and \eqref{cudiv}:
\be
\begin{aligned}\label{gdiv}
(2\mu+\lambda)&\|\sqrt{p'(\rho)^{-1}}\n\textrm{div}u\|_{L^2}^2-\frac{1}{\epsilon}\om \n \sigma\n\textrm{div}u dx\\
\le &\eta \|\n\textrm{div}u\|_{H^1}^2+C_\eta\| u\|_{H^1}^2+C\|u\|_{H^2}^2(\|u\|_{H^1}^2+\|\sigma\|_{H^2}^2 )\\
&+\frac{2\mu+\lambda}{8p'(4)}\|\n\textrm{div}u\|_{L^2}^2+C\|u_t\|_{L^2}^2.
\end{aligned}
\ee
Then we summarize \eqref{sidiv} and \eqref{gdiv} to get the lemma.

\hfill$\Box$

{\lem\label{f3}
There exist positive constants $C_5$ and $C_6$ such that
$$
\begin{aligned}
\frac{d}{dt}(\|\sqrt{\rho}&u_t\|_{L^2}^2+\|\sqrt{p'(\rho)}\sigma_t\|_{L^2}^2)+C_5\|u_t\|_{H^1}^2\\
\le &C_6\|\sigma_t\|_{H^1}^2(\|u_t\|_{L^2}^2+\|u\|_{H^1}^4+\|\sigma\|_{H^2}^2+\|\sigma_t\|_{H^1}^2)
\\
&+C_6(\epsilon^2\|\sigma_t\|_{L^2}^2+\|u\|_{H^1}^2).
\end{aligned}
$$
}
\textbf{Proof.} Note that
\begin{equation*}
\begin{split}
-\int_\Omega&(2\mu
\textrm{div} (D(u_t))+ \lambda \nabla \textrm{div}u_t)\cdot u_t dx\\
&=\om (2\mu |D(u_t)|^2 +\lambda({\rm div}u_t)^2) dx +\int_{\partial\Omega}\alpha |u_t|^2 dS\ge
\gamma_1 \|u_t\|_{H^1}^2,
\end{split}
\end{equation*}
where $\gamma_1$ is a positive constant.

We multiply \eqref{eq32} by $u_t$ and integrate to obtain
\be
\begin{aligned}\label{ut}
\frac{1}{2}&\frac{d}{dt}\|\sqrt{\rho}u_t\|_{L^2}^2+\frac{1}{\epsilon}\om p'(\rho)\nabla\sigma_t\cdot u_t dx+\gamma_1 \|u_t\|_{H^1}^2\\
&\le \eta \|u_t\|_{H^1}^2+C_\eta(\|\sigma_t\|_{H^1}^2(\|(u_t,\n \sigma)\|_{L^2}^2+\|u\|_{H^1}^4)+\|u_t\|_{H^1}^2\|u\|_{H^1}^2).
\end{aligned}
\ee

Applying $\partial_t$ to \eqref{NS3} gives
\be\label{eq416}
\sigma_{tt}+u\cdot\n \sigma_t+u_t\cdot\n \sigma+\sigma_t\textrm{div}u+\sigma\textrm{div}u_t+\frac{1}{\epsilon} \textrm{div}u_t=0.
\ee
Due to the boundary condition $u\cdot n=0$, we multiply the above equality by $p'(\rho)\sigma_t$ and integrate to get
\be
\begin{aligned}\label{st}
\frac{1}{2}\frac{d}{dt}&\|\sqrt{p'(\rho)}\sigma_t\|_{L^2}^2-\frac{1}{\epsilon}\om p'(\rho)\nabla\sigma_t\cdot u_t dx\\
\le &\eta \|u_t\|_{H^1}^2+C_\eta(\|\sigma\|_{H^2}^2\|\sigma_t\|_{L^2}^2+\|\sigma_t\|_{H^1}^4)+C\epsilon^2\|\sigma_t\|_{L^2}^2+C\|u\|_{H^1}^2.
\end{aligned}
\ee
Thus we summarize \eqref{ut} and \eqref{st} and choose $\eta$ to be sufficiently small to get the lemma.

\hfill$\Box$

\noindent Next, we should estimate the vorticity ${\rm curl}u$.
{\lem\label{f4} Let $w:={\rm curl} u$. Then
\be
\begin{aligned}\label{cu}
\frac{d}{dt}&\|\sqrt{\rho}w\|_{L^2}^2+ \mu \|{\rm{curl}} w\|_{L^2}^2\\
&\le\eta_2\|\nabla{\rm curl} w\|_{L^2}^2+ C_7\|u\|_{H^2}^2(\|u\|_{H^2}^2+\|\sigma\|_{H^2}^2)+ C_7\|u\|_{H^1}^2,
\end{aligned}
\ee
where  $\eta_2(<\mu)$ and $C_7=C_7(\eta_2)$ are  to be determined.
}\\
\textbf{Proof.} We rewrite \eqref{NS4} as
\be\label{rns4}
u_t+u\cdot\n u+\frac{1}{\epsilon}\n G( \sigma)=\frac{1}{1+\e \sigma}(\mu\triangle u+(\mu+\lambda)\n{\rm div} \textrm{u}),
\ee
where $\nabla G(\sigma)\equiv\frac{1}{\epsilon}\frac{p'(1+\e \sigma)}{1+\e \sigma}\n \sigma$ for some scalar function G.
 Applying $\lq\lq\textrm{curl}"$ to \eqref{rns4}, we obtain
\be\label{eq320}
\rho(w_t+u\cdot\n w)-\mu\triangle w=g,
\ee
where
\be\label{g}
g:=\rho h  -\epsilon\rho^{-1}\n \sigma\times(\mu\triangle u+(\mu+\lambda)\n{\rm div} \textrm{u}),
\ee
with
\be\label{h}
h:=(\partial_2u_j\partial_ju_3-\partial_3u_j\partial_ju_2,
          \partial_3u_j\partial_ju_1-\partial_1u_j\partial_ju_3,
          \partial_1u_j\partial_ju_2-\partial_2u_j\partial_ju_1)^t.
\ee
Here and in the sequel we adopt the Einstein convention about summation over repeated indices. Observing that $\Delta w=-\textrm{curl}\textrm{curl} w$, we have
\begin{equation*}
\frac{d}{dt}\|\sqrt{\rho}w\|_{L^2}^2+2\mu \|{\rm{curl}} w\|_{L^2}^2 = \om 2g\cdot w dx+2\mu\int_{\partial\Omega}({\rm{curl}}w\times n)\cdot w dS.
\end{equation*}
From \eqref{bc2}, we have
\begin{equation}\label{wbc}
\begin{aligned}
|\int_{\partial\Omega}({\rm{curl}}w\times n)\cdot w dS|
&=|\int_{\partial\Omega}(w\times n)\cdot {\rm{curl}}w dS|\\
&\le C\int_{\partial\Omega}|{\rm{curl}}w||(w\times n)\cdot\tau| dS\\
&\le \eta\|{\rm curl}w\|_{H^1}^2+C_\eta\|u\|_{H^1}^2.
\end{aligned}
\end{equation}
Thus we can easily get this lemma.

\hfill$\Box$

Next, we introduce the following two notations:
$$
\begin{aligned}
\Phi_1(t):=&\om (2\mu |D(u)|^2 +\lambda({\rm div}u)^2+|\n \sigma|^2+2 C_8(\rho |u_t|^2+p'(\rho)\sigma_t^2)\\
&+2\rho u_t\cdot u+\rho |\textrm{curl} u|^2) dx+\int_{\partial\Omega}\alpha |u|^2 dS,
\end{aligned}
$$
$$
\begin{aligned}
\Psi_1(t):=&\frac 1 2 p'(1/4)\| \sigma_t\|_{L^2}^2+(2\mu+\lambda)\|\sqrt{p'(\rho)^{-1}}\n\textrm{div}u\|_{L^2}^2\\
&+C_8C_5\|u_t\|_{H^1}^2+ \mu \|\textrm{curl} \textrm{curl} u\|_{L^2}^2,
\end{aligned}
$$
where $C_8$ is a positive constant such that $C_8>C_5^{-1}(C_2+C_3)+1$. We remark that it is important to determine the constants $C_i's$ sequentially. First, we choose $C_0,\,C_1,\,C_2,\,C_3,\,C_5$ and $C_6$ to be fixed positive constants. Next, once $\eta_1$ and $\eta_2$ are fixed, the constants $C_4$ and $C_7$ are determined.  Then it follows from Lemmas \ref{f1}, \ref{f2}, \ref{f3} and \ref{f4} that
{\lem\label{f5} Let $\epsilon_1=\min (1, \frac 1 2 \sqrt{C_8 C_6 p'(1/4)})$. Then for any $\epsilon\in (0,\epsilon_1]$,
there exists a positive constant  $C_9:=C_9(C_2,C_3,C_6,C_7,C_8)$ such that
$$
\begin{aligned}
\frac{d}{dt}  \Phi_1(t) +&\Psi_1(t)\\
\le &C_9(\|\sigma_t\|_{H^1}^2+\|u\|_{H^2}^2)(\|u_t\|_{H^1}^2+\|u\|_{H^2}^4+ \|\sigma_t\|_{H^1}^2 +\|(\sigma,u)\|_{H^2}^2)\\
+&  (C_2+2C_8C_6+C_4 +C_7)\|u\|_{H^1}^2+ \eta_1\|\n^2 {\rm{div}} u\|_{L^2}^2+\eta_2\|\n {\rm curl}{\rm curl}u\|_{L^2}^2.
\end{aligned}
$$
where $\eta_1$, $\eta_2$, $C_4(\eta_1)$ and $C_7(\eta_2)$ are positive constants to be determined later.
}
\hfill$\Box$

\subsection{The second-order estimates}
We need to estimate the  spatial and temporal derivatives of second order to close the energy estimates. The strategy is similar as that in the first-order estimates, namely, estimating the vorticity and the divergence of velocity fields respectively. However, boundary estimates are required to complete the estimates for the derivatives of highest order. We evaluate these derivatives one by one as follows.
{\lem \label{ndivu} For the positive constants $\eta_3$ and $C_{10}:=C_{10}(\eta_3)$, which are to be determined later, we have
\be\label{l7}
\begin{aligned}
(2\mu+\lambda&)\frac{d}{dt}\|\n {\rm{div}} u\|_{L^2}^2-2\frac{d}{dt}\om \rho u_t\cdot \n{\rm{div}} u dx+\|\sqrt{p'(\rho)}\n\sigma_t\|_{L^2}^2\\
\le &\eta_3(\|\n {\rm div}u_t\|_{L^2}^2+\|\n {\rm div}u\|_{H^1}^2)+C_{10}(\|u_t\|_{H^1}^2+\|\sigma_t\|_{H^1}^2\|\sigma\|_{H^2}^2\\
&+\|u\|_{H^2}^2(\epsilon^2\|\sigma_t\|_{H^1}^2\|u\|_{H^2}^2+\|u_t\|_{H^1}^2+\|\sigma\|_{H^2}^2)).
\end{aligned}
\ee
}\\
\textbf{Proof.} Note that ${\rm curl} \nabla=0$, thus
\begin{equation*}
\begin{aligned}
\om \textrm{curl} {\rm curl} u_t\cdot\n\textrm{div}u dx&= \int_{\partial\Omega}(n\times w_t)\cdot\nabla{\rm div}u dS\\
&\le C\int_{\partial\Omega} |(n\times w_t)\cdot\tau| |\nabla{\rm div}u| dS\\
&\le \eta \|\n {\rm div}u\|_{H^1}^2 + C_\eta\|u_t\|_{H^1}^2.
\end{aligned}
\end{equation*}
Multiplying both sides of \eqref{eq32} by $\n \textrm{div}u$ and integrating, we obtain
$$
\begin{aligned}
\frac{2\mu+\lambda}{2}\frac{d}{dt}&\|\n\textrm{div}u\|_{L^2}^2-\frac{d}{dt}\om \rho u_t\cdot \n\textrm{div}u dx-\frac{1}{\epsilon}\om p'(\rho)\n \sigma_t\cdot \n\textrm{div}u dx\\
\le &\eta (\|\n {\rm div}u_t\|_{L^2}^2+\|\n {\rm div}u\|_{H^1}^2)+ C_\eta\|u_t\|_{H^1}^2\\
&+C_\eta(\|u\|_{H^2}^2(\epsilon^2\|\sigma_t\|_{H^1}^2\|u\|_{H^2}^2+\|u_t\|_{H^1}^2)+\|\sigma_t\|_{H^1}^2\|\n \sigma\|_{H^1}^2).
\end{aligned}
$$
We multiply \eqref{eq37} by $p'(\rho)\n \sigma_t$ and integrate to get
\be\label{eq43}
\|\sqrt{p'(\rho)}\n\sigma_t\|_{L^2}^2+\frac{1}{\epsilon}\om p'(\rho)\n \sigma_t\cdot \n\textrm{div}u dx \le \frac 1 2\|\sqrt{p'(\rho)}\n\sigma_t\|_{L^2}^2+C \|u\|_{H^2}^2\|\sigma\|_{H^2}^2.
\ee
Combining the above two inequalities, we  get this lemma.

\hfill$\Box$


{\lem \label{n2sigma} For the positive constants $\eta_4$ and $C_{11}:=C_{11}(\eta_4)$, which are to be determined later, we have
\be\label{l8}
\begin{aligned}
\frac{d}{dt}&\|\nabla^2\sigma\|_{L^2}^2+(2\mu+\lambda)\|\sqrt{p'(\rho)^{-1}}\nabla^2{\rm div}u\|_{L^2}^2\\
&\le \eta_4\|u\|_{H^3}^2+C_{11}(\|\nabla^2 {\rm curl} u\|_{L^2}^2+\|u_t\|_{H^1}^2+\|(\sigma,u)\|_{H^2}^4+\|\sigma\|_{H^2}^2(\|u_t\|_{H^1}^2+\|u\|_{H^2}^4)).
\end{aligned}
\ee }\\
\textbf{Proof.} The following calculations are done in the form of  Einstein's convention. Applying $\partial_{ij}$ to \eqref{NS3} for $i,j=1,2,3$, where $\partial_{ij}$ denotes $\partial_{x_ix_j}$, then multiplying both sides by $\partial_{ij}\sigma$  and integrating on $\Omega$, we have
\begin{equation*}
\frac{1}{2}\frac{d}{dt}\|\partial_{ij}\sigma\|_{L^2}^2+\frac{1}{\epsilon}\om \partial_{ijk}u_k\partial_{ij}\sigma dx\le C\|\sigma\|_{H^2}^2\|u\|_{H^3}\le \eta\|u\|_{H^3}^2+C_\eta\|\sigma\|_{H^2}^4.
\end{equation*}
Next, we differentiate \eqref{NS4},   multiply the resulting equality  by $p'(\rho)^{-1}\partial_{ijk}u_k$ and integrate to obtain
\begin{equation*}
\begin{aligned}
(2\mu+\lambda&)\|\sqrt{p'(\rho)^{-1}}\partial_{ijk}u_k\|_{L^2}^2-\frac{1}{\epsilon}\om \partial_{ijk}u_k\partial_{ij}\sigma dx\\
\le & \frac{2\mu+\lambda}{2}\|\sqrt{p'(\rho)^{-1}}\partial_{ijk}u_k\|_{L^2}^2+C (\|\nabla^2 {\rm curl} u\|_{L^2}^2+ \|u_t\|_{H^1}^2\\
&+\|(\sigma,u)\|_{H^2}^4+\epsilon^2\|\sigma\|_{H^2}^2(\|u_t\|_{H^1}^2+\|u\|_{H^2}^4)).
\end{aligned}
\end{equation*}
Summarizing the above two inequalities, we show this lemma.

\hfill$\Box$

{\lem
There exists a  positive constant $\tilde{C}_{12}$ such that
\be \label{l9}
\begin{aligned}
\frac{d}{dt}(&\|{\rm div}u_t\|_{L^2}^2+\|\sqrt{\rho^{-1}p'(\rho)}\n\sigma_t\|_{L^2}^2)+(2\mu+\lambda)\|\sqrt{\rho^{-1}}\nabla{\rm div}u_t\|_{L^2}^2\\
\le &\eta_5(\|u\|_{H^3}^2+\|u_t\|_{H^2}^2)+\tilde{C}_{12} \|{\rm{curl}} w_t\|_{L^2}^2\\
&+C_{12}\|(\sigma_t,u_t)\|_{H^1}^2(\|(\sigma,u)\|_{H^2}^2+ \|u\|_{H^2}^4+\|(\sigma_t,u_t)\|_{H^1}^2),
\end{aligned}
\ee
where $\eta_5$  and ${C}_{12}:={C}_{12}(\eta_5)$ are to be chosen.
}\\
\textbf{Proof.} 
Note that with the Young inequality we have
\begin{equation*}
\om \rho^{-1}\textrm{curl} \textrm{curl} u_t\cdot\n\textrm{div}u_t dx = \eta\|\n {\rm div}u_t\|_{L^2}^2+C_\eta\|\textrm{curl} w_t\|_{L^2}^2.
\end{equation*}
By integrating the product of \eqref{eq32} and $\nabla{\rm div}u_t$, we obtain the following inequality
\be
\begin{aligned}\label{eq410}
\frac{1}{2}\frac{d}{dt}&\|\nabla{\rm div}u_t\|_{L^2}^2+(2\mu+\lambda)\|\sqrt{\rho^{-1}}\nabla{\rm div}u_t\|_{L^2}^2-\frac{1}{\epsilon}\om \rho^{-1}p'(\rho)\nabla\sigma_t\cdot\nabla{\rm div}u_t\\
\le &\frac{2\mu+\lambda}{8}\|\n {\rm div}u_t\|_{L^2}^2+C \|\textrm{curl} w_t\|_{L^2}^2
+C \|(\sigma_t,u_t)\|_{H^1}^2(\|(\sigma,u)\|_{H^2}^2+ \|u\|_{H^2}^4+\|u_t\|_{H^1}^2).
\end{aligned}
\ee
Applying $\partial_t\nabla$ to \eqref{NS3} and integrating the product of the resulting identity and $\rho^{-1}p'(\rho)\n\sigma_t$, we get
$$
\begin{aligned}
\frac{1}{2}\frac{d}{dt}&\|\sqrt{\rho^{-1}p'(\rho)}\n\sigma_t\|_{L^2}^2+\frac{1}{\epsilon}\om \rho^{-1}p'(\rho)\n\sigma_t\cdot\nabla{\rm div}u_t dx\\
=&\frac{1}{2}\om \big((\rho^{-1}p'(\rho))_t+{\rm div}(\rho^{-1}p'(\rho)u)\big)|\n \sigma_t|^2 dx\\
&-\om \rho^{-1}p'(\rho)\n \sigma_t\cdot(\n u\n\sigma_t+\n^2\sigma u_t+\n u_t\n \sigma+\n \sigma{\rm div}u_t\\
&+\sigma\n{\rm div}u_t+\n\sigma_t{\rm div}u+\sigma_t\n {\rm div}u) dx.
\end{aligned}
$$
By \eqref{NS1}, the first term on the right-hand side of the above inequality reads
\begin{equation*}
\begin{aligned}
 \frac 1 2\om (G_1(\rho)-G_1'(\rho)\rho){\rm div}u|\n \sigma_t|^2 dx
 \le \eta\|{\rm div}u\|_{H^2}^2+C_\eta\|\nabla\sigma_t\|_{L^2}^4,
\end{aligned}
\end{equation*}
where $G_1(\rho):=\rho^{-1}p'(\rho)$. Then it follows that
\be
\begin{aligned}\label{eq414}
\frac{1}{2}&\frac{d}{dt}\|\sqrt{\rho^{-1}p'(\rho)}\n\sigma_t\|_{L^2}^2+\frac{1}{\epsilon}\om \rho^{-1}p'(\rho)\n\sigma_t\cdot\nabla{\rm div}u_t dx\\
&\le \eta(\|\n u\|_{H^2}^2+\|u_t\|_{H^2}^2)+C_\eta\|\nabla\sigma_t\|_{L^2}^2(\|\sigma_t\|_{H^1}^2+\|\sigma\|_{H^2}^2).
\end{aligned}
\ee
We summarize \eqref{eq410} and \eqref{eq414} to get this lemma.

\hfill$\Box$

Next, we should derive the estimates of the vorticity $w$, which is the key of the energy estimates.
\\
{\lem \label{wh2}
There exists a  positive constant $C_{13}$ such that
\be\label{l10}
\|w\|_{H^2}^2 \le {C_{13}}(\|\Delta w\|_{L^2}^2+\|u\|_{H^2}^2).
\ee
}\\
\textbf{Proof.} Using Lemmas \ref{p2} and \ref{p3} we have
\be\label{eq524}
\|w\|_{H^2} \le C(\|{\rm{curl}}w\|_{H^1}+\|w\times n\|_{H^{\frac{3}{2}}(\partial\Omega)}+\|w\|_{H^1})
\ee
and
\be\label{eq525}
\|{\rm{curl}}w\|_{H^1} \le C(\|\Delta w\|_{L^2}+\|{\rm{curl}}w\cdot n\|_{H^{\frac{1}{2}}(\partial\Omega)}+\|{\rm{curl}}w\|_{L^2}).
\ee
From  \eqref{bc2} and the trace theorem, we obtain
\be\label{eq526}
\|w\times n\|_{H^{\frac{3}{2}}(\partial\Omega)} \le C\|u\|_{H^{\frac{3}{2}}(\partial\Omega)} \le C\|u\|_{H^2}.
\ee
We construct the local coordinates by the isothermal coordinates $\lambda(\psi,\varphi)$ to derive an estimate near the boundary (see \cite{JO} for instance), where $\lambda(\psi,\varphi)$ satisfies
$$\lambda_\psi\cdot\lambda_\psi>0, \; \lambda_\varphi\cdot\lambda_\varphi>0 \; {\rm and}\;  \lambda_\psi\cdot\lambda_\varphi=0.$$
We cover the boundary $\partial\Omega$ by a finite number of bounded open sets $W^k\subset\mathbb{R}^3, k=1,2,...,L,$ such that for any $x \in W^k\cap\Omega,$
$$x = \lambda^k(\psi,\varphi)+rn(\lambda^k(\psi,\varphi)) = \Lambda^k(\psi,\varphi,r),$$
where $\lambda^k(\psi,\varphi)$ is the isothermal coordinate and  $n$ is the unit outer normal to $\partial\Omega$.
For simplicity, in what follows we will omit the superscript $k$ in each $W^k$. Then we construct the orthonormal system corresponding to the local  coordinates by
$$e_1 = \frac{\lambda_\psi}{|\lambda_\psi|},\; e_2 = \frac{\lambda_\varphi}{|\lambda_\varphi|},\; e_3 = n(\lambda) = e_1\times e_2.$$
By a straightforward calculation, we see that $J \in C^2$ and
$$\begin{aligned}
J &= det Jac\Lambda = (\Lambda_\psi\times\Lambda_\varphi)\cdot e_3\\
&=|\lambda_\psi||\lambda_\varphi|+r(|\lambda_\psi|n_\varphi\cdot e_2+|\lambda_\varphi|n_\psi\cdot e_1)+r^2[(n_\psi\cdot e_1)(n_\varphi\cdot e_2)-(n_\psi\cdot e_2)(n_\varphi\cdot e_1)] > 0,
\end{aligned}
$$
for sufficiently small $r>0$.
Obviously, $Jac(\Lambda^{-1}) = (Jac\Lambda)^{-1}$. Moreover, we can easily derive the following relations (see also \cite{V}):
$$[\nabla(\Lambda^{-1})^1]\circ\Lambda = \frac{1}{J}(\Lambda_\psi\times e_3),$$
$$[\nabla(\Lambda^{-1})^2]\circ\Lambda = \frac{1}{J}(e_3 \times \Lambda_\varphi),$$
$$[\nabla(\Lambda^{-1})^3]\circ\Lambda = \frac{1}{J}(\Lambda_\varphi\times\Lambda_\psi),$$
where the notation $'\circ'$ stands for the composition of operators. Set $y:=(y_1,y_2,y_3):=(\psi,\varphi,r)$, $a_{ij} = ((Jac\Lambda)^{-1})_{ij}$. Then $n = (a_{31},a_{32},a_{33})$, the tangential directions $\tau_i = (a_{i1},a_{i2},a_{i3})(i=1,2)$, and
 \be
a_{ij}a_{3j}=0,\quad{\rm for}\;i=1,2.
\ee
 Then we denote by $D_i$ the partial derivative with respect to $y_i$ in local coordinates. To be precise, $D_3$ is the normal derivative and $D_i$ for $i=1,2$ are the tangential derivatives in the original coordinates. Moreover, we have
$$\partial_{x_j}=a_{kj}D_k.$$
 Next, we denote the vorticity near the boundary as $\widetilde{w}:=(\widetilde{w}_1,\widetilde{w}_2,\widetilde{w}_3)^t:=w(t,\Lambda(y))$. By direct calculations we get
$$\begin{aligned}
{\rm{curl}}w\cdot n =& (a_{k2}D_k\widetilde{w}_3-a_{k3}D_k\widetilde{w}_2,a_{k3}D_k\widetilde{w}_1-a_{k1}D_k\widetilde{w}_3,a_{k1}D_k\widetilde{w}_2-a_{k2}D_k\widetilde{w}_1)\cdot(a_{31},a_{32},a_{33})\\
=& [(a_{32}a_{13}-a_{33}a_{12})D_1+(a_{32}a_{23}-a_{33}a_{22})D_2]\widetilde{w}_1\\
&+[(a_{33}a_{11}-a_{31}a_{13})D_1+(a_{33}a_{21}-a_{31}a_{23})D_2]\widetilde{w}_2\\
&+[(a_{31}a_{12}-a_{32}a_{11})D_1+(a_{31}a_{22}-a_{32}a_{21})D_2]\widetilde{w}_3\\
=& \sum_{i=1}^2(n\times\tau_i)\cdot D_i\widetilde{w}\\
=& \sum_{i=1}^2(D_i((n\times\tau_i)\cdot \widetilde{w})-D_i(n\times\tau_i)\cdot \widetilde{w})\\
=&\sum_{i=1}^2(D_i((n\times \widetilde{w})\cdot \tau_i)-D_i(n\times\tau_i)\cdot \widetilde{w}).
\end{aligned}
$$
Thus, with \eqref{eq526} we obtain that
\be\label{eq527}
\begin{aligned}
\|{\rm{curl}}w\cdot n\|_{H^{\frac{1}{2}}(\partial\Omega)} &\leq \sum_{i=1}^2\|D_i((n\times w)\cdot \tau_i)\|_{H^{\frac{1}{2}}(\partial\Omega)}+C\sum_{i=1}^3\| w_j\|_{H^{\frac{1}{2}}(\partial\Omega)}\\
&\le  C\|u\|_{H^2}.
\end{aligned}
\ee
With $\Delta w=-\textrm{curl}\textrm{curl} w$ it follows from \eqref{eq524}, \eqref{eq525}, \eqref{eq526} and \eqref{eq527} that
\be\label{curlu}
\|w\|_{H^2} \le {C}(\|\Delta w\|_{L^2}+\|u\|_{H^2}).
\ee
\hfill$\Box$
{\lem
There exists a positive constant $C_{14}$  such that
\be\label{deltaw}
\begin{aligned}
\mu&\frac{d}{dt}\|{\rm curl} w\|_{L^2}^2+\frac{\mu^2}{20}\|\Delta w\|_{L^2}^2+\|\sqrt {\rho}w_t\|_{L^2}^2\\
&\le \eta_6\|u\|_{H^2}^2+C_{14}(\|u\|_{H^2}^4+\|\sigma\|_{H^2}^2\|u\|_{H^3}^2)+C_{15}\|u_t\|_{H^1}^2,
\end{aligned}
\ee
where $\eta_6$ and $C_{15}:=C_{15}(\eta_6)$ are to be chosen.
}\\

\noindent\textbf{Proof.} Note that
$$ \om \Delta ww_t dx = -\frac{1}{2}\frac{d}{dt}\om |{\rm curl} w|^2 dx + \int_{\partial\Omega}(n\times {\rm curl} w)\cdot w_t dS $$
and
\begin{eqnarray*}
\int_{\partial\Omega}(n\times {\rm curl} w)\cdot w_t dS & =&\int_{\partial\Omega} {\rm curl} w\cdot
(n\times w_t) dS\\
&\le& \eta\|{\rm curl}w\|_{L^2(\p\Omega)}^2+ C_\eta \|(n\times w_t)\cdot\tau\|_{L^2(\p\Omega)}^2\\
&\le &  \eta\|{\rm curl}w\|_{H^1}^2 + C_\eta\|u_t\|_{H^1}^2
\end{eqnarray*}
invoking of Lemma \ref{lem25}.
Multiplying \eqref{eq320} by $w_t-\delta\Delta w$, where $\delta$ is a positive constant to be chosen, and integrating, we get
\begin{equation*}
\begin{aligned}
\frac{\mu}{2}\frac{d}{dt}&\|{\rm curl} w\|_{L^2}^2+\mu\delta\|\Delta w\|_{L^2}^2+\|\sqrt {\rho}w_t\|_{L^2}^2\\
\le &\frac{1}{4}\|\sqrt {\rho}w_t\|_{L^2}^2+C(\|g\|_{L^2}^2+\|u\|_{H^2}^2\|\n w\|_{L^2}^2)+\frac{1}{4}\|\sqrt {\rho}w_t\|_{L^2}^2+\delta^2\|\rho\|_{L^\infty}\|\Delta w\|_{L^2}^2\\
&+\delta^2\|\Delta w\|_{L^2}^2+C(\|u\|_{H^2}^2\|\n w\|_{L^2}^2+\|g\|_{L^2}^2)+\eta\|{\rm curl} w\|_{H^1}^2+C_\eta\|u_t\|_{H^1}^2.
\end{aligned}
\end{equation*}
In virtue of \eqref{g},\eqref{h} and \eqref{curlu}, we get this lemma by choosing $\eta$ small enough and $\delta=\frac{\mu}{10}$.

\hfill$\Box$

{\lem\label{wt}
There exists a positive constant  $C_{16}$ such that
\be\label{l12}
\begin{aligned}
\frac{d}{dt}\|&\sqrt {\rho}w_t\|_{L^2}^2+\mu \|{\rm{curl}} w_t\|_{L^2}^2\\
\le &\eta_7(\|u_{tt}\|_{H^1}^2+\|u_t\|_{H^2}^2) + C_{16}\|u_t\|_{H^1}^2 \\ + & C_{17}(\|\sigma_t\|_{H^1}^2(\|u_t\|_{H^1}^2+\|u\|_{H^2}^4+\|u\|_{H^3}^2+\|\sigma\|_{H^2}^2\|u\|_{H^2}^2)
+\|u_t\|_{H^2}^2\|(u,\sigma)\|_{H^2}^2),
\end{aligned}
\ee
where $\eta_6$ and $C_{17}:=C_{17}(\eta_6)$ are to be chosen.
}\\
\textbf{Proof.} From \eqref{eq320}, we have
\be\label{eq522}
\rho (w_{tt} +u\cdot \nabla w_t)-\mu\triangle w_t=\widetilde{g},
\ee
where $\widetilde{g}:=g_t-\epsilon\sigma_t(w_t+u\cdot\n w)-\rho u_t\cdot\n w$ with
$$
\begin{aligned}
|g_t| \le &C(\epsilon|\sigma_t||\n u|^2+|\n u_t||\n u|+\epsilon^2|\sigma_t||\n \sigma||\nabla^2 u|\\
&+\epsilon(|\sigma_t||\n^2 u|+|\n \sigma||\n^2 u_t|)).
\end{aligned}
$$
Multiplying \eqref{eq522} by $w_t$ and integrating over $\Omega$, we have
$$
\frac{d}{dt}\|\sqrt {\rho}w_t\|_{L^2}^2+2\mu \|{\rm curl} w_t\|_{L^2}^2 = \om 2\widetilde{g}\cdot w_t dx + 2\mu \int_{\partial\Omega}{\rm curl} w_t\cdot (w_t\times n) dS.
$$
Similar as \eqref{wbc} and \eqref{eq525}, one has
$$
\begin{aligned}
\int_{\partial\Omega}{\rm curl} w_t\cdot (w_t\times n) dS
\le \eta\|{\rm{curl}} w_t\|_{H^1}^2 + C_\eta\|u_t\|_{H^1}^2.
\end{aligned}
$$
and
$$
\|{\rm{curl}}w_t\|_{H^1} \le C(\|\Delta w_t\|_{L^2}+\|{\rm{curl}}w_t\cdot n\|_{H^{\frac{1}{2}}(\partial\Omega)}+\|{\rm{curl}}w_t\|_{L^2}).
$$
From \eqref{eq522} again we get
$$
\mu \|\Delta w_t\|_{L^2}^2 \le C(\|u_{tt}\|_{H^1}^2+\|\widetilde{g}\|_{L^2}^2+ \|u\|_{H^2}^2\|u_t\|_{H^1}^2).
$$
Moreover, similar as \eqref{eq527}, we  can derive that
$$ \|{\rm{curl}}w_t\cdot n\|_{H^{\frac{1}{2}}(\partial\Omega)} \le  C\|u_t\|_{H^2}.$$
Collecting all the above estimates, this lemma is shown.

\hfill$\Box$

To close the energy estimates, we estimate $\sigma_{tt}$ and $u_{tt}$ in the following two lemmas.
{\lem\label{s7}
There exists a positive constant $C_{18}$  such that
\be\label{l13}
\begin{aligned}
\frac{d}{dt} \|\e( &\sqrt{p'(\rho) }\sigma_{tt},\sqrt{\rho}u_{tt})\|_{L^2}^2+C_{18}\|\e u_{tt}\|_{H^1}^2\\
 &\le   \eta_8 \|({\rm div}u,u_t)\|_{H^2}^2+C_{19}[\|\sigma_t\|_{H^1}^4(1+ \|\sigma\|_{H^2}^2)+\|u_t\|_{H^2}^2(\|\sigma_t\|_{H^1}^2\|u\|_{H^2}^2+\|u_t\|_{H^1}^2) \\
&+\|\e(\sigma_{tt},u_{tt})\|_{L^2}^2(\|(\sigma_t,u_t)\|_{H^1}^2 + \|(\sigma,u)\|_{H^2}^2+\|u\|_{H^2}^4+\|\e\sigma_{tt}\|_{L^2}^2)],
\end{aligned}
\ee
where $\eta_8$ and $C_{19}:=C_{19}(\eta_8)$ are to be determined.
}\\
\textbf{Proof.} We differentiate \eqref{NS3} twice with respect to $t$, multiply the resulting equality   by $\e^2p'(\rho)\sigma_{tt}$ and then integrate over $\Omega$. Finally we get
\begin{equation}\label{sigmatt}
\begin{aligned}
\frac{1}{2}&\frac{d}{dt}\|\epsilon\sqrt{p'(\rho)}\sigma_{tt}\|_{L^2}^2-\e\om p'(\rho)\n \sigma_{tt}\cdot u_{tt} dx\\
&\le \eta \|({\rm div}u,u_t)\|_{H^2}^2+ \delta\|\e u_{tt}\|_{H^1}^2 +C_{\eta,\delta}(\|\e\sigma_{tt}\|_{L^2}^2+\|\sigma_t\|_{H^1}^2+\|\sigma\|_{H^2}^2)\|\e\sigma_{tt}\|_{L^2}^2.
\end{aligned}
\end{equation}
On the other hand, we get the following equality by differentiating \eqref{eq32} in temporal variable:
\be
\begin{aligned}\label{eq528}
\rho(u_{ttt}+u\cdot\n u_{tt})+\frac{1}{\epsilon}p'(\rho)\n \sigma_{tt}-2\mu{\rm div}(D(u_{tt}))- \lambda \n{\textrm{div}}u_{tt}=f,
\end{aligned}
\ee
where
$$
\begin{aligned}
-f=&2p''\sigma_t\n \sigma_t+(\epsilon p'''\sigma_t^2+p''\sigma_{tt})\n\sigma+\epsilon\sigma_{tt}u_t+2\epsilon\sigma_tu_{tt}\\
&+(\rho u_{tt}+\epsilon\sigma_{tt}u+2\epsilon\sigma_tu_t)\cdot\n u+2(\epsilon\sigma_tu+\rho u_t)\cdot\n u_t.
\end{aligned}
$$
Multiplying \eqref{eq528} by $\e^2u_{tt}$, then integrating on $\Omega$, we get
\begin{equation}\label{utt}
\begin{aligned}
\frac{1}{2}\frac{d}{dt} \om & \rho|\e u_{tt}|^2 dx+\e\om p'(\rho)\n \sigma_{tt}\cdot u_{tt} dx\\
&+\om (2\mu |\e D(u_{tt})|^2+\lambda (\e{\rm{div}} u_{tt})^2) dx +\int_{\partial\Omega}\alpha |\e u_{tt}|^2 dS\\
&\le \delta\|\e u_{tt}\|_{H^1}^2+ C_\delta [\|\e(\sigma_{tt},u_{tt})\|_{L^2}^2(\|(\sigma_t,u_t)\|_{H^1}^2 + \|(\sigma,u)\|_{H^2}^2+\|u\|_{H^2}^4) \\
&+\|\sigma_t\|_{H^1}^4(1+ \|\sigma\|_{H^2}^2)+\|u_t\|_{H^2}^2(\|\sigma_t\|_{H^1}^2\|u\|_{H^2}^2+\|u_t\|_{H^1}^2)].
\end{aligned}
\end{equation}
Note that all the above calculations can be verified rigorously by regularization arguments. Summarizing \eqref{sigmatt} and \eqref{utt} and selecting $\eta$ small enough, we can prove this lemma.

\hfill$\Box$

{\lem\label{s8}
There exists a positive constant $C_{20}$  such that
\be\label{l14}
\begin{aligned}
\frac{d}{dt}\om \e^2\rho u_{tt}&\cdot u_t dx +\frac{\e^2}{2}\frac{d}{dt}(\om (2\mu |D(u_t)|^2 +\lambda({\rm div}u_t)^2) dx +\int_{\partial\Omega}\alpha |u_t|^2 dS)+C_{20}\|\e\sigma_{tt}\|_{L^2}^2\\
\le &\eta_9\|u_t\|_{H^1}^2+5\|u_{tt}\|_{H^1}^2+C_{21}\e(\|\sigma_t\|_{H^1}^4(1+\|\sigma\|_{H^1}^2) +\|\sigma_t\|_{H^1}^2\|u\|_{H^2}^2\\
&+ \|u_t\|_{H^1}^2(\|(\sigma_t,u_t)\|_{H^1}^2 + \|(\sigma,u)\|_{H^2}^2 +\|\sigma_t\|_{H^1}^2\|u\|_{H^2}^2 +\|u\|_{H^2}^4 )),
\end{aligned}
\ee
where $\eta_9$ and $C_{21}:=C_{21}(\eta_9)$ are to be determined.
}\\
\textbf{Proof.} We integrate $\e^2p'(\rho)\sigma_{tt}$ times \eqref{eq416} to get
\begin{equation*}
\begin{aligned}
\|\e\sqrt{p'(\rho)} \sigma_{tt}\|_{L^2}^2& +2\e\om p'(\rho){\textrm{div}}u_t\sigma_{tt} dx\\
&\le   C \e(\|u\|_{H^2}^2\|\sigma_t\|_{H^1}^2+\|u_t\|_{H^1}^2\|\sigma\|_{H^2}^2).
\end{aligned}
\end{equation*}
Multiplying \eqref{eq528} by $\e^2u_t$ and integrate over $\Omega$, we obtain
\begin{equation*}
\begin{aligned}
\frac{d}{dt}\om \e^2 &\rho u_{tt}\cdot u_t dx  + \e\om p'(\rho)\n \sigma_{tt}\cdot u_t dx + \om \e^2 p''(\rho)\n \sigma\cdot u_t\sigma_{tt} dx\\
+ & \frac{\e^2}{2}\frac{d}{dt}(\om  (2\mu |D(u_t)|^2 +\lambda({\rm div}u_t)^2) dx +\int_{\partial\Omega}\alpha |u_t|^2 dS)\\
\le &  \om \rho |\e u_{tt}|^2 dx  + \epsilon^3 |\om \sigma_t u_{tt}\cdot u_t dx|
+|\om \e^2 p''(\rho)\n \sigma\cdot u_t\sigma_{tt} dx| + |\om \e^2 f \cdot u_t dx|\\
\le &  5\|\e u_{tt}\|_{H^1}^2 +\eta\|(\e\sigma_{tt}, u_t) \|_{H^1}^2 + C_\eta\e(\|\sigma_t\|_{H^1}^4(1+\|\sigma\|_{H^1}^2) \\
+ & \|u_t\|_{H^1}^2(\|(\sigma_t,u_t)\|_{H^1}^2 + \|(\sigma,u)\|_{H^2}^2 +\|\sigma_t\|_{H^1}^2\|u\|_{H^2}^2 +\|u\|_{H^2}^4 )).
\end{aligned}
\end{equation*}
Thus summarizing the above inequalities we get this lemma.
\hfill$\Box$

{\lem\label{lem314}
There exists a positive constant $C_{22}$ such that
\be\label{eq539}
\|\sigma\|_{H^2}^2 \le C_{22}\epsilon(\|u_t\|_{H^1}^2+\|u\|_{H^2}^4+\|u\|_{H^3}^2)(1+\|\sigma\|_{H^1}^2).
\ee
}\\
\textbf{Proof.} From \eqref{NS3} and \eqref{NS5}, we deduce that
$$\frac{d}{dt}\om \sigma dx = - \int_{\partial\Omega} (\sigma+\frac{1}{\epsilon})u\cdot n dS = 0.$$
By the assumption $\om \sigma_0 dx = 0$, we have $\om \sigma dx = 0$. Then this lemma follows from \eqref{NS4} and the Poincar\'e inequality
\begin{equation*}
\begin{aligned}
\|\sigma\|_{H^2} &\le C\|\n \sigma\|_{H^1}.
\end{aligned}
\end{equation*}

\hfill$\Box$

We introduce the following notations:
\begin{equation*}
\Phi(t):=C_{23}\Phi_0(t)+C_{24}\Phi_1(t)+\Phi_2(t),
\end{equation*}
\begin{equation*}
\Psi(t):=C_{23}\Psi_0(t)+C_{24}\Psi_1(t)+\Psi_2(t),
\end{equation*}
where
$$\begin{aligned}
\Phi_2(t):=&(2\mu+\lambda)\|\n {\rm{div}} u\|_{L^2}^2-2\om \rho u_t\cdot \n{\rm{div}} u dx\\
&+\|\nabla^2\sigma\|_{L^2}^2+\|{\rm div}u_t\|_{L^2}^2+\|\sqrt{\rho^{-1}p'(\rho)}\n\sigma_t\|_{L^2}^2\\
&+\frac{\widetilde{C}_{12}+2}{\mu}\|\sqrt {\rho}{\rm{curl}}u_t\|_{L^2}^2+\frac{20(C_{11}+2)C_{13}}{\mu}\|{\rm{curl}} {\rm{curl}}u\|_{L^2}^2\\
&+\om \e^2\rho u_{tt}\cdot u_t dx+K\|\e(\sqrt{p'(\rho)}\sigma_{tt},\sqrt{\rho}u_{tt})\|_{L^2}^2\\
&+\frac{\e^2}{2}(\om (2\mu |D(u_t)|^2 +\lambda({\rm div}u_t)^2) dx +\int_{\partial\Omega}\alpha |u_t|^2 dS)
\end{aligned}$$
and
$$\begin{aligned}
\Psi_2(t):=&(2\mu+\lambda)(\|\sqrt{p'(\rho)^{-1}}\nabla^2{\rm div}u\|_{L^2}^2+\|\sqrt{\rho^{-1}}\nabla{\rm div}u_t\|_{L^2}^2)\\
&+\|\sqrt{p'(\rho)}\n\sigma_t\|_{L^2}^2+2\|{\rm{curl}} {\rm{curl}}u_t\|_{L^2}^2\\
&+\frac{20(C_{11}+2)C_{13}}{\mu^2}\|\sqrt {\rho}{\rm{curl}}u_t\|_{L^2}^2+2\|{\rm curl} u\|_{H^2}^2\\
&+C_{20}\|\e\sigma_{tt}\|_{L^2}^2+2\|\e u_{tt}\|_{H^1}^2+\|\sigma\|_{H^2}^2.
\end{aligned}$$
 Here we choose $K$ satisfying $K>7C_{18}^{-1}$ and the constants $C_{23}$ and $C_{24}$ are to be determined.
We summarize \eqref{l7}, \eqref{l8}, \eqref{l9}, $(C_{11}+2)\times$\eqref{l10}, $20\mu^{-2} (C_{11}+ 2)C_{13}\times$\eqref{deltaw}, $\mu^{-1}(\tilde{C}_{12}+2)\times$\eqref{l12}, $K\times$\eqref{l13}, \eqref{l14} and \eqref{eq539},   with the aid of
Lemmas \ref{p2} and \ref{p3}, 
and then choose $\eta_i(i=3,4,5,6,7,8,9)$ and $\epsilon$ small enough in the resulting inequality.
 Then we get the estimates of highest order as follows:
{\lem\label{s6}
There exist positive constants $\epsilon_2\in(0,\epsilon_1]$ and $C_{26}$ such that
\be\label{eq531}
\begin{aligned}
\frac{d}{dt}\Phi&_2(t)+\frac{3}{4}\Psi_2(t)\\
\le& (C_{10}+C_{11}+\frac{20(C_{11}+2)C_{13}}{\mu^2}C_{15}+\frac{\widetilde{C}_{12}+2}{\mu}C_{16}+1)\|u_t\|_{H^1}^2\\
+ &((C_{11}+2)C_{13}+1)\|u\|_{H^2}^2
 +C_{26}\Psi(t)(\Phi(t)+\Phi(t)^2).
\end{aligned}
\ee
}


Next we  redefine the constant $C_8$ to close the energy estimates:
$$
\begin{aligned}
C_8  & \geq  4C_5^{-1}(C_2+C_3+C_5(\frac{32}{2\mu+\lambda}+\frac{4}{K})+C_{10}+C_{11}\\
&  +\frac{20(C_{11}+2)C_{13}}{\mu^2}C_{15}+\frac{\widetilde{C}_{12}+2}{\mu}C_{16}+2).
 \end{aligned}
$$
Let
\begin{equation*}
C_{24}\ge 4\overline{C}(1+(C_{11}+2)C_{13})/\min(2\mu+\lambda,\mu),
\end{equation*}
for some constant $\overline C$ depending on the best constants in Lemmas \ref{p2} and \ref{p3}, next choose $\eta_1$ and $\eta_2$ small enough in Lemma \ref{f5} and then set
\begin{equation*}
C_{25}:=C_{24}(C_2+2C_8C_6+C_4 +C_7+\overline{C}(1+(C_{11}+2)C_{13})).
\end{equation*}
Thus we have
$$
\begin{aligned}
\frac{d}{dt} (C_{24}\Phi_1 & +\Phi_2) + \frac{3}{4}(C_{24}\Psi_1+\Psi_2)\\  &\le  C\Psi(t)(\Phi(t)+\Phi(t)^2)+C_{25}\|u\|_{H^1}^2.
\end{aligned}
$$
Finally, we choose $C_{23}$ such that
\begin{equation*}
C_{23}C_0 \geq C_{25},
\end{equation*}
and next $\e$ small enough in Lemma \ref{z1}.
Then there exist positive constants  $\epsilon_3\in (0,\epsilon_2]$ and $C \in [1,+\infty)$ such that
$$
\frac{d}{dt} \Phi(t) +\frac{1}{2}\Psi(t) \le \frac{1}{2}C \Psi(t)(\Phi(t)+\Phi(t)^2), \quad \forall\, 0\le t\le T, \quad 0 \le \epsilon \le \epsilon_3.
$$
That is,
\be\label{eq542}
\frac{d}{dt}\Phi(t) \le -\frac{1}{2}\Psi(t)(1-C\Psi(t)(\Phi(t)+\Phi(t)^2)), \quad\forall \,0\le t\le T, \quad 0 \le \epsilon \le \epsilon_3.
\ee
Then we can obtain the following lemma, which can be shown exactly in the same way as in {\cite{V,OU1}}. Thus the details are omitted.
{\lem {\rm(}Uniform Estimates{\rm)} \label{lem58}
Let $\Omega\subset\mathbb{R}^3$ be a simply connected, bounded domain with smooth boundary $\partial\Omega$. Let $(u,\sigma)$ be a solution to \eqref{NS3}-\eqref{NS6} in $\Omega\times(0,T)$ with $\frac{1}{4} \le 1+\epsilon\sigma \le 4$, $\forall (x,t)\in\Omega\times(0,T)$, $\epsilon \in (0,\epsilon_3]$. Suppose that
$$
\Phi(0) \le \frac{\beta}{\Theta}, \quad \beta\in(0,\frac{1}{2}],
$$
for some constant $\Theta>0.$ Then we have
$$
\Phi(t) \le \frac{\beta}{\Theta},  \quad \forall\, t\in[0,T].
$$
}
\hfill$\Box$

\noindent{\bf Acknowledgement.} This work was partially supported by NSFC under grant 11001021, the China Postdoctoral Science Foundation under grant 201003077, and the Fundamental Research Funds for the Central Universities.





\bibliographystyle{elsarticle-num}







\end{document}